\documentclass{article}
\usepackage{graphicx, latexsym,amsmath,amssymb}
\def\margine#1{\strut\vadjust{\kern-\strutdepth\specialstar{#1}}}

\def\strutdepth{\dp\strutbox}
\def\specialstar#1{\vtop to \strutdepth{\baselineskip\strutdepth
    \vss\llap{#1\quad}\null}}
\def\Proof{\smallbreak\par\noindent{\sc Proof.~\hskip1pt}}

\let\qq=\qquad

\def\texIT#1{\indent\llap{#1\enspace}\ignorespaces}
 \def\IT#1 {\smallbreak\hangindent\parindent\texIT{#1}
    \ifdim\lastskip=\smallskipamount\removelastskip\penalty55\
    \smallskip\fi}

\newcount\nequa
\def\e(#1){\ifmmode\eqlabel{#1}\eqno(\the\nequa)
          \ifbozze\rlap{\qq #1}\fi\else(\eqref{#1})\fi}
\def\ee(#1)#2{\ifmmode\eqlabel{#1}\eqno(\the\nequa{\rm#2})
          \ifbozze\rlap{\qq #1}\fi\else(\eqref{#1}#2)\fi}
\def\z(#1){\eqno(\the\nequa{\rm #1})}
\def\eqref#1{\reference{#1}{e}}
\def\eqlabel#1{\global\advance\nequa by1%
               \scrivi{#1}{e}{EQUAZIONE}{\nequa}}
\def\zz(#1){(\nequa{\rm #1})
           \ifbozze\rlap{\qq\etichetta}\fi}
\def\eqal(#1)#2{\ifbozze\xdef\etichetta{#1}\fi\eqlabel{#1}\eqalignno{#2}}

\newif\ifbozze

\def\reference#1#2{%
    \seindefinito{@#2@#1}%
     \else%
       \csname @#2@#1\endcsname%
     \fi}

\def\seindefinito#1{\expandafter\ifx\csname #1\endcsname\relax}

\def\scrivi#1#2#3#4{%
     \xdef\ilnumero{\the#4}%
     \expandafter\xdef\csname  @#2@#1\endcsname{\ilnumero}%
    }

\newcount\nconjecture

\def\Conjecture(#1){\vskip-0.2\baselineskip\relax\goodbreak%
\advance\nconjecture by 1%
\vspace{2mm}
\ET{\bf Conjecture \the\nconjecture}%
       \ifbozze\margine{$\Diamond\;$\scriptsize#1}\fi
    \scrivi{#1}{CN}{CONJECTURE}{\nconjecture}
     \begin{sl}}
\def\fineconjecture{\end{sl}%
\subsection*{}\vskip-1.2\baselineskip\relax\goodbreak%
}

\def\endconjecture{\end{sl}\vspace{2mm}\goodbreak\ET}

\def\conjecture(#1){Conjecture~\hskip-2pt\reference{#1}{CN}}

\newcount\ncorollary

\def\Corollary(#1){\vskip-0.2\baselineskip\relax\goodbreak%
\advance\ncorollary by 1%
\vspace{2mm}
\ET{\bf Corollary \the\ncorollary}%
       \ifbozze\margine{$\Diamond\;$\scriptsize#1}\fi
    \scrivi{#1}{CO}{COROLLARY}{\ncorollary}
     \begin{sl}}
\def\finecorollary{\end{sl}%
\subsection*{}\vskip-1.2\baselineskip\relax\goodbreak%
}

\def\endcorollary{\end{sl}\vspace{2mm}\goodbreak\ET}

\def\corollary(#1){Corollary~\hskip-2pt\reference{#1}{CO}}

\newcount\nlemma

\def\Lemma(#1){\vskip-0.2\baselineskip\relax\goodbreak%
\advance\nlemma by 1%
\vspace{2mm}
\ET{\bf Lemma \the\nlemma}%
       \ifbozze\margine{$\Diamond\;$\scriptsize#1}\fi
    \scrivi{#1}{LE}{LEMMA}{\nlemma}
     \begin{sl}}
\def\finelemma{\end{sl}%
\subsection*{}\vskip-1.2\baselineskip\relax\goodbreak%
}

\def\endlemma{\end{sl}\vspace{2mm}\goodbreak\ET}

\def\lemma(#1){Lemma~\hskip-2pt\reference{#1}{LE}}

\newcount\ndefinition

\def\Definition(#1){\advance\ndefinition by 1
\vspace{2mm}
\ET{\bf Definition \the\ndefinition}%
       \ifbozze\margine{$\Diamond\;$\scriptsize#1}\fi
    \scrivi{#1}{D}{DEFINITION}{\ndefinition}
    \begin{sl} }

\def\ET{\smallskip\par\noindent}

\def\enddefinition{\end{sl}\vspace{2mm}\goodbreak\ET}

\def\definition(#1){Definition~\hskip-0pt\reference{#1}{D}}

\newcount\nteorema

\def\Theorem(#1){\advance\nteorema by 1
\vspace{2mm}
\ET{\bf Theorem \the\nteorema}%
       \ifbozze\margine{$\Diamond\;$\scriptsize#1}\fi
    \scrivi{#1}{T}{TEOREMA}{\nteorema}
    \begin{sl} }

\def\ET{\smallskip\par\noindent}

\def\endtheorem{\end{sl}\vspace{2mm}\goodbreak\ET}

\def\theorem(#1){Theorem~\hskip-0pt\reference{#1}{T}}

\newcount\nproposition

\def\Proposition(#1){\advance\nproposition by 1
\vspace{2mm}
\ET{\bf Proposition \the\nproposition}%
       \ifbozze\margine{$\Diamond\;$\scriptsize#1}\fi
    \scrivi{#1}{P}{PROPOSITION}{\nproposition}%
    \begin{sl} }
 
\def\ET{\smallskip\par\noindent}

\def\endproposition{\end{sl}\vspace{0mm}\goodbreak\ET}

\def\proposition(#1){Proposition~\hskip-0pt\reference{#1}{P}}

\newdimen\pIR
\def\R{{\rm I\kern\pIR R}}

\newdimen\pIN
\def\N{{\rm I\kern\pIN N}}

\newdimen\pIC
\def\C{{\rm I\kern\pIC C}}

\newdimen\pIS
\def\S{{\rm I\kern\pIS S}}

\def\Proof{\noindent{\em Proof.}\enspace}

\def\endproof{\hfill $\Box$ \vspace{3mm}}

\begin{document}

\title{\bf Excessive $[l,m]$-factorizations}
\author{David Cariolaro\\
Department of Mathematical Sciences\\
Xi'an Jiaotong-Liverpool University\\
Suzhou Higher Education Town\\
Suzhou, Jiangsu\\
215123 China\\
\\
E-mail: \texttt{david.cariolaro@xjtlu.edu.cn}\\
\\
Giuseppe Mazzuoccolo \thanks{Research supported by a fellowship from the
    European Project ``INdAM fellowships in mathematics and/or
    applications for experienced researchers cofunded by Marie Curie
    actions''} \\Laboratoire G-SCOP (Grenoble-INP, CNRS)\\
    Grenoble\\ France\\
\\
E-mail: \texttt{mazzuoccolo@unimore.it}}

\maketitle

\begin{abstract} 
\noindent Given two positive integers $l$ and $m$, with $l\leq m$, an \textit{$[l,m]$-covering} of a graph $G$ is a set $\mathcal{M}$ of matchings of $G$ whose union is the edge set of $G$ and such that $l\leq |M|\leq m$ for every $M\in {\mathcal{M}}$.

\noindent An $[l,m]$-covering $\mathcal{M}$ of $G$ is an \textit{excessive $[l,m]$-factorization} of $G$ if the cardinality of $\mathcal{M}$ is as small as possible.
The number of matchings in an excessive $[l,m]$-factorization of $G$ (or $\infty,$ if $G$ does not admit an excessive $[l,m]$-factorization) is a graph parameter called the \textit{excessive $[l,m]$-index} of $G$ and denoted by $\chi'_{[l,m]}(G)$. In this paper we study such parameter. Our main result is a general formula for the excessive $[l,m]$-index of a graph $G$ in terms of other graph parameters. Furthermore, we give a polynomial time algorithm which computes $\chi'_{[l,m]}(G)$ and outputs an excessive $[l,m]$-factorization of $G$, whenever the latter exists.

\end{abstract}

\vspace{3mm}

\bozzefalse
\section{Introduction}
The classical concept of graph factorization as the decomposition of the edge set of a graph into (pairwise isomorphic) factors is a very general concept which has received a substantial amount of attention in the literature. One limitation of the use of such concept is that it is normally applicable only to specific classes of graphs, such as complete graphs, or $k$-factorizable graphs, etc. In 2004 one extension of the concept of $1$-factorization, called \textit{excessive factorization,} which is applicable to a wider class of graphs, has been proposed \cite{DAV} (see also \cite{BER}).
Informally speaking, an excessive factorization of a graph $G$ is a minimum set of (not necessarily edge-disjoint) $1$-factors of $G$ whose union is the edge set of $G$. Thus, in order for a graph to admit an excessive factorization, it is not necessary that it is $1$-factorizable (or even regular), and hence, using this new concept, one can develop and apply the results of the corresponding theory to a much wider class of graphs.
Of course one may observe that there are limitations also in the concept of excessive factorization, in what it applies only to graphs having $1$-factors and, more precisely, having $1$-factors containing any prescribed edge of the graph. It is therefore desirable to study extensions of this concept by replacing the term ``$1$-factor" by something more general. However, if we replace the term ``$1$-factor" by ``arbitrary matching" what we obtain is essentially the concept of \textit{edge colouring}, which has been studied since the nineteenth century and is therefore not a new concept.
An intermediate possibility is to replace the term ``1-factor" by ``matching of fixed size $m$", and this idea was pursued by Cariolaro and Fu in \cite{CF}, where the corresponding concept was called ``excessive $[m]$-factorization". 

More precisely an \textit{excessive $[m]$-factorization} of a graph $G$ is a set ${\mathcal{M}}$ of matchings of $G$ such that  
\begin{description}
\item (i) $\bigcup_{M\in {\mathcal{M}}}M=E(G)$;
\item (ii) $|M|=m$ for every $M\in {\mathcal{M}}$;
\item (iii) subject to (i) and (ii), $|{\mathcal{M}}|$ is minimum.
\end{description}

A set ${\mathcal{M}}$ of matchings of $G$ satisfying conditions (i) and (ii) above, but not necessarily  (iii), is called an \textit{$[m]$-covering} of $G$.
A graph which admits an $[m]$-covering is said to be \textit{$[m]$-coverable}.
It is obvious that a graph $G$ admits an excessive $[m]$-factorization if and only it is $[m]$-coverable, which is the case if and only if every edge $e$ of $G$ belongs to a matching of size $m$ (or, equivalently, at least $m$) of $G$. Such condition can be verified in polynomial time thanks to a famous theorem of Edmonds \cite{EDM}. The number of matchings in an excessive $[m]$-factorization (or $\infty$, if $G$ does not admit an excessive $[m]$-factorization) is a graph parameter which is denoted in \cite{CF} by $\chi'_{[m]}(G)$ and called the \textit{excessive $[m]$-index} of $G$. 

The theory of excessive factorizations is still in its infancy, but a number of papers have already been written on the topic (see e.g. \cite{BER,BM,MULT,COM3,BR,MAZ1}) and connections with some important combinatorial problems such as the Berge-Fulkerson Conjecture have already been noticed \cite{MAZZ}.

Whilst finding an excessive factorization in general is an NP-hard problem \cite{BER}, it was recently established by Cariolaro and Rizzi \cite{CR} that, for a fixed value of $m$, there exists a polynomial time algorith which, given as input a graph $G$, outputs the excessive $[m]$-index $\chi'_{[m]}(G)$ as well as an excessive $[m]$-factorization. 

The purpose of this paper is to introduce a generalization of the concept of excessive $[m]$-factorization, as follows. Let $l,m$ be two positive integers, where $l\leq m$. An \textit{excessive $[l,m]$-factorization} of $G$ is a set $\mathcal{M}$ of matchings of $G$ such that
\begin{description}
\item (i) $\bigcup_{M\in {\mathcal{M}}}M=E(G)$;
\item (ii) $l\leq |M|\leq m$ for every $M\in {\mathcal{M}}$;
\item (iii) subject to (i) and (ii), $|{\mathcal{M}}|$ is minimum.
\end{description}
A set ${\mathcal{M}}$ of matchings of $G$ satisfying conditions (i) and (ii) above, but not necessarily  (iii), is called an \textit{$[l,m]$-covering} of $G$. A graph is said to be \textit{$[l,m]$-coverable} if it admits an $[l,m]$-covering. For notational convenience, a matching $M$ satisfying $l\leq |M|\leq m$ will be called an \textit{$[l,m]$-matching} and, in the case $l=m,$ it will simply be called an \textit{$[m]$-matching}.

Similarly to the case of excessive $[m]$-factorizations, we define \textit{excessive $[l,m]$-index} of the graph $G$, denoted by $\chi'_{[l,m]}(G),$ as the cardinality of an excessive $[l,m]$-factorization of $G$ if $G$ admits an excessive $[l,m]$-factorization, and $\infty$ otherwise.

Notice that, when $l=m$, the concepts of excessive $[l,m]$-factorization and excessive $[l,m]$-index coincide, respectively, with the concepts of excessive $[m]$-factorization and excessive $[m]$-index.

Our main result (Theorem 3) will be a general formula for the excessive $[l,m]$-index of a graph $G$ expressed in terms of the chromatic index of $G$ and the excessive $[k]$-index of $G$, for some particular values of the integer $k$. A natural question is whether, for a fixed value of the integers $l$ and $m$, where $l\leq m,$ there exists a polynomial time algorithm which, given a graph $G$, computes $\chi'_{[l,m]}(G)$ and outputs an excessive $[l,m]$-factorization of $G$. We prove in the last section that the answer to this question is affirmative.

\section{Preliminary results and definitions}
An edge colouring of a multigraph $G$ is a map $\varphi: E(G)\rightarrow {\mathcal{C}}$, where $\mathcal{C}$ is a set (called the set of \textit{colours}) and $\varphi$ has the property of mapping adjacent edges into distinct colours.
When $|{\mathcal{C}}|=k,$ $\varphi$ is called a \textit{$k$-edge colouring}. A \textit{colour class} of $\varphi$ is a set of edges of the form $\varphi^{-1}(\{\alpha\}),$ where $\alpha$ is a colour. 
The \textit{chromatic index} of $G$, denoted by $\chi'(G),$ is the minimum integer $k$ such that $G$ has a $k$-edge colouring.

A $k$-edge colouring $\varphi$ is called an \textit{equalized $k$-edge colouring} if, for every colour class $C$ of $\varphi,$ we have
$$
\lfloor |E(G)|/k \rfloor \leq |C|\leq \lceil |E(G)|/k \rceil.
$$
The following result, obtained independently by McDiarmid \cite{MC} and de Werra \cite{DW}, will be used often in the sequel.

\Lemma(DW) Let $G$ be a multigraph and suppose $G$ has a $k$-edge colouring. Then $G$ admits an equalized $k$-edge colouring. Furthermore an equalized $k$-edge colouring can be found in time $O(|V||E|)$.
\endlemma

We shall also need the following lemma of Cariolaro and Fu \cite[Theorem 6]{CF}.

\Lemma(CF) Let $G$ be a graph and let $m$ be an integer such that $|E(G)|/m\geq \chi'(G)$. Then $\chi'_{[m]}(G)=\lceil |E(G)|/m \rceil$.
\endlemma

Let $l,m$ be two integers, with $l\leq m,$ and let $\mathcal{M}$ be an $[l,m]$-covering of $G$. The \textit{multigraph $\tilde{G}$ induced by $\mathcal{M}$} is the multigraph with the same vertex set as $G$, where two distinct vertices $u,v$ are joined by as many edges in $\tilde{G}$ as there are matchings in $\mathcal{M}$ containing the edge $uv$. Similarly, if $\mathcal{H}$ is any multigraph whose underlying simple graph is $G$, and if $\varphi$ is a $k$-edge colouring of $\mathcal{H}$, with colour classes $C_1,C_2,\ldots, C_k$, the \textit{covering of $G$ induced by $\varphi$} is the covering
${\mathcal{M}}=\{M_1,M_2,\ldots, M_k\},$ where $M_i$ is the matching of $G$ defined by
$$
M_i=\{uv\in E(G)\;|\: \mbox{ there exists $e\in C_i$ such that $e$ joins $u$ and $v$ in ${\mathcal{H}}$.}\}
$$

Henceforward, whenever it is not specified, the symbols $l$ and $m$ will denote two positive integers, satisfying $l\leq m$.
 
We have the following.

\Proposition(1) The graph $G$ admits an excessive $[l,m]$-factorization if and only if it admits an excessive $[l]$-factorization.
\endproposition
\Proof Every $[l]$-covering of $G$ is also an $[l,m]$-covering of $G$. Hence the existence of an excessive $[l]$-factorization implies the existence of an $[l,m]$-factorization. Conversely, if $G$ admits an excessive $[l,m]$-factorization, then, in particular, every edge of $G$ belongs to a matching of size at least $l$, and hence $G$ admits an excessive $[l]$-factorization.\endproof

\Proposition(1b) For every positive integers $l,l',m,m',$ with $l'\leq l\leq m\leq m'$ and every graph $G$, we have $\chi'_{[l',m']}(G)\leq \chi'_{[l,m]}(G)$.
\endproposition
\Proof Obvious since every $[l,m]$-covering of $G$ is an $[l',m']$-covering of $G$. \endproof

The following proposition generalizes \cite[Proposition 1]{CF}.

\Proposition(3) The following conditions are equivalent for any graph $G$.
\begin{description}
\item (i) $\chi'_{[l,m]}(G)\leq k;$
\item (ii) $G$ has a $k$-edge colouring $\varphi$ such that each colour class of $\varphi$ is contained in an $[l,m]$-matching of $G$;
\item (iii) $G$ is the underlying simple graph of a multigraph $\tilde{G}$ which is $k$-edge colourable and whose colour classes are $[l,m]$-matchings of $\tilde{G}$.
\end{description}
\endproposition
\Proof Assume (i). Let ${\mathcal{M}}=\{M_1,M_2,\ldots ,M_k\}$ be an $[l,m]$-covering of $G$, where, if necessary, we allow the  same matching to appear more than once in $\mathcal{M}$. Define a function $\varphi: E(G)\rightarrow \{1,2,\ldots, k\}$ by $\varphi(e)=\min_{1\leq i\leq k}\{i\;|\; e\in M_i\}.$
It is straightforward to verify that $\varphi$ is an edge colouring of $G$ whose colour classes can each be extended to an $[l,m]$-matching of $G$. This shows that (i) implies (ii).

Assume now (ii). Let $\varphi$ be a $k$-edge colouring whose colour classes are contained in an $[l,m]$-matching of $G$. Let $N_1,N_2,\ldots,N_k$ be the colour classes of $\varphi$.
By assumption, for every $N_i$ there is an $[l,m]$-matching $M_i$ of $G$ containing $N_i$. Thus $\{M_1,M_2,\ldots,M_k\}$ is an $[l,m]$-covering of $G$, whence (i) follows.

Assume now (i), and let ${\mathcal{M}}=\{M_1,M_2,\ldots,M_k\}$ be an $[l,m]$-covering of $G$, as above. Let $\tilde{G}$ be the multigraph induced by $\mathcal{M}$. By construction, $\tilde{G}$ has a $k$-edge colouring whose colour classes are $[l,m]$-matchings of $\tilde{G}$, hence (iii) follows.

Conversely, if $\tilde{G}$ has a $k$-edge colouring $\psi$ whose colour classes $\{C_1,C_2,\ldots,C_k\}$ are $[l,m]$-matchings of $G$, it suffices to consider the $[l,m]$-covering of $G$ induced by $\psi$.
Clearly $\mathcal{M}$ is an $[l,m]$-covering of $G$, whence (i) follows. Thus (i),(ii) and (iii) are equivalent.
\endproof

\Proposition(4) $\chi'_{[l,m]}(G)\geq \max\{\chi'(G), \lceil \frac{|E(G)|}{m} \rceil\}.$
\endproposition
\Proof Let $k=\chi'_{[l,m]}(G)$. We can assume that $k$ is finite. By \proposition(3) (ii), $G$ has a $k$-edge colouring, hence $k\geq \chi'(G)$.
Let ${\mathcal{M}}=\{M_1,M_2,\ldots, M_k\}$ be an excessive $[l,m]$-factorization of $G$. Since each matching in $\mathcal{M}$ has size at most $m$, we have
$$
|E(G)|=|\bigcup_{i=1}^kM_i|\leq km,
$$
hence $k\geq |E(G)|/m, $ and since $k$ is an integer, we obtain $k\geq \lceil |E(G)|/m \rceil$. This concludes the proof. \endproof

\section{Proof of the main result}
In this section we assume that the integers $l$ and $m$ satisfy the inequality $l<m$, unless stated otherwise.
We have the following.

\Lemma(A) $\chi'_{[l,m]}(G)=\min_{l\leq i<m}\chi'_{[i,i+1]}(G)$.
\endlemma
\Proof By \proposition(1b), we have
$$
\chi'_{[l,m]}(G)\leq \min_{l\leq i<m}\chi'_{[i,i+1]}(G).
$$
We now prove the reverse inequality. In doing so, we can clearly assume that $\chi'_{[l,m]}(G)<\infty$. Let $\mathcal{M}$ be an excessive $[l,m]$-factorization of $G$, and assume $|{\mathcal{M}}|=k$. Let $\tilde{G}$ be the multigraph induced by $\mathcal{M}$. Notice that
$$
kl\leq |E(\tilde{G})|\leq km, \e(Y)
$$
and, by construction, $\tilde{G}$ is $k$-edge colourable. By \lemma(DW), $\tilde{G}$ has an equalized $k$-edge colouring $\varphi$.
In such colouring, every colour class has size $\lfloor \frac{|E(\tilde{G})|}{k} \rfloor$ or $\lceil \frac{|E(\tilde{G})|}{k} \rceil.$
By \e(Y), 
$$
\lfloor \frac{|E(\tilde{G})|}{k} \rfloor \geq l
$$
and 
$$
\lceil \frac{|E(\tilde{G})|}{k}  \rceil \leq m.
$$
Hence, letting 
$$
i=\lfloor \frac{|E(\tilde{G})|}{k}  \rfloor,
$$
the edge colouring $\varphi$ of $\tilde{G}$ induces an $[i,i+1]$-covering of $G$ of cardinality $k$, thus proving
$$
\chi'_{[i,i+1]}(G)\leq k.
$$
This terminates the proof. \endproof 

\Lemma(B) If the integer $i$ satisfies $i\geq \frac{|E(G)|}{\chi'(G)},$ then $\chi'_{[i,i+1]}(G)=\chi'_{[i]}(G)$.
\endlemma
\Proof By \proposition(1b), $\chi'_{[i,i+1]}(G)\leq \chi'_{[i]}(G)$. We prove the reverse inequality.
We can clearly assume that $\chi'_{[i,i+1]}(G)=k<\infty$.
Let $\mathcal{M}$ be an excessive $[i,i+1]$-factorization of $G$ and let $\lambda_i$ (respectively, $\lambda_{i+1}$) be the number of $[i]$-matchings (respectively, $[i+1]$-matchings) in $\mathcal{M}$.
Let $\tilde{G}$ be the multigraph induced by $\mathcal{M}$. By \proposition(3) (iii), $\tilde{G}$ is $k$-edge colourable.
Notice that
$$
|E(\tilde{G})|=i\lambda_i+(i+1)\lambda_{i+1}=i(\lambda_i+\lambda_{i+1})+\lambda_{i+1}=i\chi'_{[i,i+1]}(G)+\lambda_{i+1}
$$
$$
\geq i\chi'(G)+\lambda_{i+1}\geq |E(G)|+\lambda_{i+1},
$$
where in the proof of the last inequality we have used our assumption that $i\geq |E(G)|/\chi'(G)$. 
Thus, in particular, we can delete $\lambda_{i+1}$ edges from $\tilde{G}$ and still obtain a multigraph $\tilde{H}$ which has $G$ as its underlying simple graph. By definition, $\tilde{H}$ contains 
$$i(\lambda_i+\lambda_{i+1})=i\chi'_{[i,i+1]}(G)$$
edges and is $k$-edge colourable (since $\tilde{G}$ is).
Let $\varphi$ be an equalized $k$-edge colouring of $\tilde{H}$ (which exists by \lemma(DW)).
Then $\varphi$ induces a covering of $G$ with $k$ matchings of size $i$, thus proving
$$
\chi'_{[i]}(G)\leq k.
$$
\endproof

\Lemma(C) If the integer $i$ satisfies $i\geq |E(G)|/\chi'(G)$, then $\chi'_{[i+1]}(G)\geq \chi'_{[i]}(G).$
\endlemma
\Proof Without loss of generality, we may assume that $\chi'_{[i+1]}(G)=k<\infty$.
Let $\mathcal{M}$ be an excessive $[i+1]$-factorization of $G$, and let $\tilde{G}$ be the multigraph induced by $\mathcal{M}$.
We have
$$
|E(\tilde{G})|=k(i+1)=ki+k\geq \chi'(G)i+k\geq |E(G)|+k,
$$
where in the last inequality we have used the assumption.
Hence we may delete $k$ edges from $\tilde{G}$ and still obtain a multigraph $\tilde{H}$ whose underlying simple graph is $G$.
Notice that
$$
|E(\tilde{H})|=ki
$$
and $\tilde{H}$ is $k$-edge colourable, since $\tilde{G}$ is $k$-edge colourable (by \proposition(3) (iii)).
Let $\varphi$ be an equalized $k$-edge colouring of $\tilde{H}$. Clearly $\varphi$ induces an $[i]$-covering of $G$ with $k$ matchings, therefore proving that $\chi'_{[i]}(G)\leq k$.
This terminates the proof.\endproof

\Lemma(D) If the integer $l$ satisfies $l\geq \frac{|E(G)}{\chi'(G)},$ then $\chi'_{[l,m]}(G)=\chi'_{[l]}(G)$.
\endlemma
\Proof Using \lemma(A), \lemma(B) and \lemma(C), we have
$$
\chi'_{[l,m]}(G)=\min_{l\leq i<m}\chi'_{[i,i+1]}(G)=\min_{l\leq i<m}\chi'_{[i]}(G)=\chi'_{[l]}(G),
$$
as desired.
\endproof

We are now in a position to prove our main result.

\Theorem(X) For every pair of positive integers $l,m$ with $l\leq m$, and any graph $G$, we have
$$\chi'_{[l,m]}(G)= \begin{cases} 
\lceil \frac{|E(G)|}{m} \rceil & \text{ if } \frac{|E(G)|}{\chi'(G)}\geq m \\
\\
\chi'(G) & \text{ if } l \le \frac{|E(G)|}{\chi'(G)} \leq m\\
\\
\chi'_{[l]}(G) & \text{ if } \frac{|E(G)|}{\chi'(G)} \leq l \\
\end{cases} $$
\endtheorem
\Proof First observe that the result holds for $l=m$ by \lemma(CF). We now assume $l<m$. Suppose first that 
$$
|E(G)|/\chi'(G)\geq m.
$$ 
It follows from \lemma(CF) that 
$$\chi'_{[m]}(G)=\lceil |E(G)|/m \rceil.
$$
By \proposition(4) and \proposition(1b), we have
$$
\lceil |E(G)|/m \rceil \leq \chi'_{[l,m]}(G)\leq \chi'_{[m]}(G)=\lceil |E(G)|/m \rceil,
$$
hence we have $$\chi'_{[l,m]}(G)=\lceil |E(G)|/m \rceil.$$

Suppose now that 
$$
l\leq |E(G)|/\chi'(G) \leq m.
$$

Let $k=\chi'(G)$ and let $\varphi$ be an equalized $k$-edge colouring of $G$, with colour classes $C_1,C_2,\ldots,C_k$.
Notice that, for every $i=1,2,\ldots,k,$ we have
$$
l\leq \lfloor |E(G)|/\chi'(G) \rfloor \leq |C_i| \leq \lceil |E(G)|/\chi'(G) \rceil \leq m,
$$
hence $\{C_1,C_2,\ldots, C_k\}$ is an $[l,m]$-covering of $G$, which implies $
k\geq \chi'_{[l,m]}(G).$ Hence $\chi'(G)\geq \chi'_{[l,m]}(G)$. 
The reverse inequality  follows from \proposition(4).

Suppose now 
$$|E(G)|/\chi'(G)\leq l.$$
By \lemma(D), we have
$$
\chi'_{[l,m]}(G)=\chi'_{[l]}(G),
$$
as desired. This concludes the proof. \endproof

\section{Extremal cases}
Consider the special case $l=1.$ In this case \theorem(X) reduces to the following.

\Corollary(R) $\chi'_{[1,m]}(G)=\max\{\chi'(G),\lceil |E(G)|/m \rceil\}.$
\endcorollary

Clearly an excessive $[1,m]$-factorization $\mathcal{M}$ of a graph $G$ is just a minimum set of matchings of size at most $m$ whose union is $E(G)$.
Thus (since we are only interested in minimum coverings) there is clearly no loss of generality in assuming that the matchings in $\mathcal{M}$ are disjoint, and hence that $\mathcal{M}$ is an edge colouring whose colour classes have size at most $m$. Such colouring was called an \textit{optimal $m$-bounded edge colouring} in a recent paper of Rizzi and the first author \cite{RC}, where \textit{inter alia} it was shown that, for a fixed value of the integer $m$, an optimal $m$-bounded edge colouring of any graph $G$ (and hence the parameter $\chi'_{[1,m]}(G)$) can be computed in polynomial time (see Theorem $5$ in Section $7$).

A further extremal case is obtained by considering the case $m=\infty$, i.e. we consider coverings with matchings of size at least $l$ but with no prescribed upper bound on their size. Attention to this case was prompted to us by Richard Brualdi (oral communication with the first author at the 7th Shanghai Conference on Combinatorics in 2011).
The corresponding factorization is called an \textit{excessive $[l,\infty]$-factorization} and the corresponding parameter is called an \textit{excessive $[l,\infty]$-index} and denoted by $\chi'_{[l,\infty]}(G)$.

We notice that, in general, the problem of the computation of this parameter is NP-hard, since it is easily seen that $\chi'_{[1,\infty]}(G)=\chi'(G)$, and it is well known that computing $\chi'(G)$ is NP-hard \cite{HOL}. 
The following result follows from \theorem(X).

\Corollary(Y) For every integer $l$ and any graph $G$, we have
$$\chi'_{[l,\infty]}(G)= \begin{cases} 
\chi'(G) & \text{ if }  \frac{|E(G)|}{\chi'} \geq l\\
\\
\chi'_{[l]}(G) & \text{ if } \frac{|E(G)|}{\chi'} \leq l \\
\end{cases} $$
\endcorollary

\section{Compatibility}
\proposition(4) shows that
$$
\chi'_{[l,m]}(G)\geq \max\{\chi'(G),\lceil |E(G)|/m\rceil \}.
$$
Thus, in particular,
$$
\chi'_{[m]}(G)\geq \max\{\chi'(G),\lceil |E(G)|/m\rceil \}.
$$
Graphs for which the above inequality holds as an equality were called \textit{$[m]$-compatible} in \cite{CF}.
It was proved in \cite{CF} that, for every graph $G$, there exists an integer $com(G)$, called \textit{compatibility index}, such that
$$
\mbox{$G$ is $[m]$-compatible if and only if $1\leq m\leq com(G).$}
$$
Generalizing this notion, we say that $G$ is \textit{$[l,m]$-compatible} if 
$$
\chi'_{[l,m]}(G)=\max\{\chi'(G),\lceil |E(G)|/m\rceil.
$$
This definition naturally suggests the following question: for a fixed graph $G$ and integer $m$, for which values of $l$ is $G$ $[l,m]$-compatible?

It follows from \corollary(R) that $G$ is always $[1,m]$-compatible.

Suppose now that $G$ is $[l,m]$-compatible and $l>1$. Then
$$
\chi'_{[l-1,m]}(G)\leq \chi'_{[l,m]}(G)=\max\{\chi'(G), \lceil \frac{|E(G)|}{m}\rceil\},
$$
and hence, using \proposition(4), we see that $G$ is $[l-1,m]$-compatible. Thus, for every $m$ there is an integer $f_G(m)$ such that $G$ is $[l,m]$-compatible if and only if $1\leq l\leq f_G(m)$. In particular $f_G(m)=m$ if and only if $G$ is $[m]$-compatible.
We call the function $f_G: {\mathbb{Z}}^+\rightarrow {\mathbb{Z}}^+$ the \textit{compatibility function} of $G$.

\begin{figure}[h]
\begin{center}
\includegraphics[width=60mm]{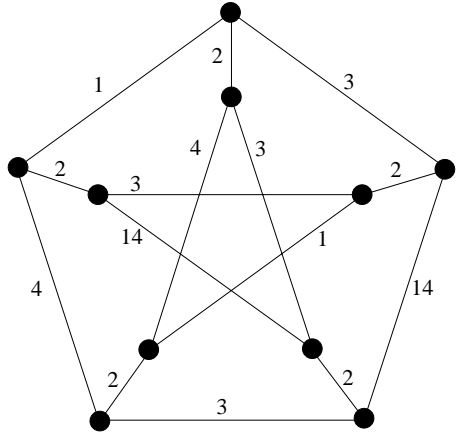}
\caption{An excessive $[4,5]$-factorization of the Petersen graph}
\label{pet45}
\end{center}
\end{figure}

For example, if $P$ is the Petersen graph, since it is known \cite{CF} that $com(P)=4,$ it follows that $f_P(m)=m$ for every $m\leq 4$ and $f_P(5)<5$.
In Fig. \ref{pet45} we show a $[4,5]$-covering of $P$ consisting of $4$ matchings, hence necessarily an excessive $[4,5]$-factorization, thereby proving that $P$ is $[4,5]$-compatible and that $f_P(5)=4$.

We now prove that the function $f_G$ is always nondecreasing.

\Theorem(Z) Let $G$ be a graph. Then the compatibility function $f_G$ is nondecreasing.
\endtheorem
\Proof It will cleary suffice to prove that, if $G$ is $[l,m]$-compatible, for two integers $l$ and $m$, then it is $[l,m+1]$-compatible.
Suppose that $G$ is $[l,m]$-compatible. Let ${\mathcal{M}}=\{M_1,M_2\ldots, M_k\}$ be an excessive $[l,m]$-factorization. Notice that, by definition,
$$
k=\max\{\chi'(G), \lceil |E(G)|/m \rceil\}.
$$
Let 
$$
k'=\max\{\chi'(G), \lceil |E(G)|/(m+1) \rceil\}.
$$
Notice that $k'\leq k$. If $k'=k,$ then necessarily $\mathcal{M}$ is an excessive $[l,m+1]$-factorization.
Therefore we can assume that $k'<k.$ We now divide the proof in two cases.

\noindent \textbf{Case 1: $k'=\lceil \frac{|E(G)|}{m+1} \rceil$.}

In this case
$$
\chi'(G)\leq k'=\lceil |E(G)|/(m+1)\rceil < \lceil |E(G)|/m \rceil=k. \e(U1)
$$
In particular, $G$ is $k'$-edge colourable. Since
$$
|E(G)|/(m+1)\leq k',
$$
we have
$$
\lceil |E(G)|/k'\rceil \leq m+1.\e(U2)
$$
Since 
$$
|E(G)|/m >k',
$$
and since $l\leq m, $ we have 
$$
|E(G)|/l>k',
$$
and hence
$$
|E(G)|/k'>l,
$$
so that
$$
\lfloor |E(G)|/k' \rfloor \geq l. \e(U3)
$$
By \e(U2) and \e(U3), an equalized $k'$-edge colouring of $G$ is an $[l,m+1]$-covering, and hence necessarily an excessive $[l,m+1]$-factorization.
Thus 
$$
\chi'_{[l,m+1]}(G)=k',
$$
and hence $G$ is $[l,m+1]$-compatible. 

\noindent \textbf{Case 2: $k'=\chi'(G)>\lceil \frac{|E(G)|}{m+1}\rceil $.}

Since $k>k',$ we have 
$$k=\lceil \frac{|E(G)|}{m} \rceil >k'=\chi'(G)>\lceil \frac{|E(G)}{m+1}\rceil.$$

We need to prove that $G$ is $[l,m+1]$-compatible. Notice that
$$
\frac{|E(G)|}{l}\geq \frac{|E(G)|}{m}>\chi'(G)>\frac{|E(G)}{m+1}.
$$
Let $\varphi$ be an equalized $\chi'(G)$-edge colouring. Then every colour class $C$ satisfies
$$
l\leq m\leq \lfloor \frac{|E(G)}{\chi'(G)}\rfloor \leq |C|\leq \lceil \frac{|E(G)|}{\chi'(G)}\rceil \leq m+1,
$$
hence $\varphi$ is an excessive $[l,m]$-factorization, and we conclude that $G$ is $[l,m+1]$-compatible.\endproof

\section{Coherence}
We have the following.

\Proposition(2) $\chi'_{[l,m]}(G)\leq \min_{l\leq i\leq m}\chi'_{[i]}(G)$.
\endproposition
\Proof 
Let $i$ be an integer, with $l\leq i\leq m$. Without loss of generality we may assume that $\chi'_{[i]}(G)$ is finite. Let $\mathcal{M}$ be an excessive $[i]$-factorization of $G$. Then $\mathcal{M}$ is also an $[l,m]$-covering of $G$, implying that $\chi'_{[l,m]}(G)\leq \chi'_{[i]}(G)$. By the arbitrariety of $i$, the assertion is proved.
\endproof

\begin{figure}[h]
\begin{center}
\includegraphics[width=50mm]{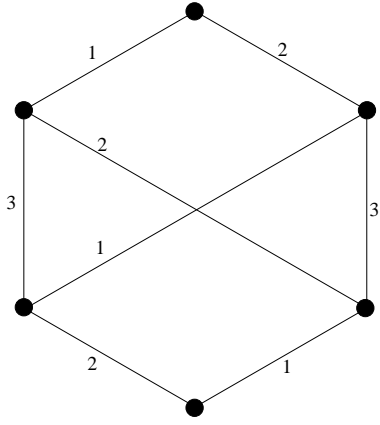}
\caption{A graph which is not $[2,3]$-coherent. As shown in the figure, $\chi'_{[2,3]}=3.$ It is easy to see that $\chi'_{[2]}(G)=\chi'_{[3]}(G)=4$.}
\label{inco}
\end{center}
\end{figure}

A graph for which the inequality expressed by \proposition(2) holds as an equality will be called $[l,m]$-\textit{coherent}. Notice that every graph is $[m,m]$-coherent by definition.
An example of a graph $G$ and two integers $l,m$ such that $G$ is not $[l,m]$-coherent is shown in Fig. \ref{inco}.

The following theorem gives a characterization of the graphs which are not $[l,m]$-coherent.

\Theorem(COH) A graph $G$ is not $[l,m]$-coherent if and only if $l < \frac{|E(G)|}{\chi'(G)} < m$ and $\chi'_{[k]}(G)>\chi'(G),$ where $k=\lceil \frac{|E(G)|}{\chi'(G)} \rceil.$ 	
\endtheorem
\Proof Assume that $G$ is not $[l,m]$-coherent. Then clearly $l<m$. By \theorem(X), we have 
$$
\frac{|E(G)|}{\chi'(G)}>l. \e(HA)
$$
Assume 
$$
\frac{|E(G)|}{\chi'(G)}\geq m.\e(H1)
$$
Then, by \lemma(CF) and \theorem(X), we have

$$
\chi'_{[l,m]}(G)= \lceil \frac{|E(G)|}{m} \rceil = \chi'_ {[m]}(G),
$$
hence $G$ is $[l,m]$-coherent, a contradiction. Therefore \e(H1) is false, and we have
$$
\frac{|E(G)|}{\chi'(G)}< m.\e(HB)
$$
By \e(HA), \e(HB) and \theorem(X) we then have
$$
\chi'_{[l,m]}(G)=\chi'(G).
$$
Since $G$ is not $[l,m]$-coherent,
$$
\chi'_{[i]}(G)>\chi'_{[l,m]}(G)=\chi'(G) \mbox{ for every $i, l\leq i\leq m$.}
$$
In particular, letting $k=\lceil \frac{|E(G)|}{\chi'(G)}\rceil,$
by \e(HA) and \e(HB) we have $l\leq k\leq m$, and by \e(HC) we have
$$
\chi'_{[k]}(G)>\chi'(G),
$$
as desired.

Suppose now that $G$ is $[l,m]$-coherent and assume 
$$
l<\frac{|E(G)|}{\chi'(G)}<m
$$
and 
$$
\chi'_{[k]}(G)>\chi'(G),
$$
where $k=\lceil \frac{|E(G)|}{\chi'(G)}\rceil.$

By \theorem(X), we have
$$
\chi'_{[l,m]}(G)=\chi'(G).
$$
Moreover, by \lemma(C), if $i \geq k$, then 
$$\chi'_{[i+1]}(G)\geq \chi'_{[k]}(G)> \chi'(G).$$
On the other hand, if $i<\frac{|E(G)|}{\chi'(G)},$ then by \lemma(CF)
$$
\chi'_{[i]}(G)=\lceil \frac{|E(G)|}{i}\rceil > \chi'(G).
$$
Thus $\chi'_{[i]}(G)>\chi'(G)$ for every $i, l\leq i\leq m,$ and hence $G$ is not $[l,m]$-coherent, a contradiction.
This contradiction concludes the proof.\endproof

The example of Fig. \ref{inco} shows that there are $[l,m]$-compatible graphs which are not $[l,m]$-coherent.
On the other hand, there are graphs which are $[l,m]$-coherent and not $[l,m]$-compatible. For example any graph which is not $[m]$-compatible (e.g. the Petersen graph for $m=5$) is nonetheless $[m,m]$-coherent.

\section{Complexity}
We shall now prove that, for any fixed positive integers $l,m$, with $l\leq m,$ there exists a polynomial time algorithm that, given a graph $G$, outputs $\chi'_{[l,m]}(G)$ and, if $\chi'_{[l,m]}(G)<\infty,$ also outputs an excessive $[l,m]$-factorization of $G$.

First notice that, using \corollary(R) and \proposition(4), we can restate our \theorem(X) as follows.

\Theorem(XX) For any pair of positive integers $l,m$, with $l\leq m,$ and any graph $G$, we have
$$\chi'_{[l,m]}(G)= \begin{cases} 
\chi'_{[1,m]} & \text{ if } \frac{|E(G)|}{\chi'(G)}\geq l \\
\\
\chi'_{[l]}(G) & \text{ if } \frac{|E(G)|}{\chi'(G)} \leq l \\
\end{cases} $$
\endtheorem

We shall use the following two results of Rizzi and the first author, which we have already mentioned but which, for convenience, we state below.

\Theorem(RC) \cite{RC} Let $m$ be a fixed positive integer. Then there exists a polynomial time algorithm which, given a graph $G$, outputs $\chi'_{[1,m]}(G)$ as well as an excessive $[1,m]$-factorization of $G$.
\endtheorem

Notice that we can always assume that the excessive $[1,m]$-factorization obtained as a result of \theorem(RC) is an edge colouring, whose colour classes are all of size at most $m$ (optimal $m$-bounded edge colouring).

\Theorem(CR) \cite{CR} Let $m$ be a fixed positive integer. Then there exists a polynomial time algorithm which, given a graph $G$, outputs $\chi'_{[m]}(G)$ and, if $\chi'_{[m]}(G)<\infty,$ also outputs an excessive $[m]$-factorization of $G$.
\endtheorem

Our  algorithm, which we name $EXC(G,l,m)$, is outlined below.

\noindent \textbf{ALGORITHM $EXC(G,l,m)$}
\begin{enumerate}
\item INPUT $G$.
\item Compute (using \theorem(RC)) $\chi'_{[1,l]}(G)$, $\chi'_{[1,m]}(G)$ and an $m$-bounded edge colouring  $\varphi$ of $G$.
\item IF $\chi'_{[1,m]}(G)<\chi'_{[1,l]}(G),$ then transform $\varphi$ in an equalized edge colouring $\varphi'$ using Mc Diarmid and de Werra's algorithm (\lemma(DW)).
\item RETURN $\chi'_{[1,m]}$ and $\varphi'$.
\item ELSE compute $\chi'_{[l]}(G)$ and, if $\chi'_{[l]}(G)<\infty$, compute an excessive $[l]$-factorization $\mathcal{M}$ of $G$ using \theorem(CR).
\item RETURN $\chi'_{[l]}(G)$ and (if $\chi'_{[l]}(G)<\infty$) $\mathcal{M}$.
\end{enumerate}

We shall now prove that Algorithm $EXC(G,l,m)$ is correct.

\Theorem(CM) For any fixed positive integers $l,m$, with $l\leq m,$ Algorithm $EXC(G,l,m)$ computes in polynomial time $\chi'_{[l,m]}(G)$ and, if $\chi'_{[l,m]}(G)<\infty,$ outputs an excessive $[l,m]$-factorization of $G$. 
\endtheorem
\Proof 
By \theorem(XX) we have that $\chi'_{[l,m]}(G)$ equals either $\chi'_{[1,m]}(G)$ or $\chi'_{[l]}(G)$. Since an $[l]$-covering of $G$ is an $[1,m]$-covering of $G$, the relation $\chi'_{[1,m]}(G) \leq \chi'_{[l]}(G)$ holds.
Suppose now $\chi'_{[1,m]}(G) < \chi'_{[1,l]}(G)$. Let $\varphi$ be an optimal $m$-bounded edge colouring of $G$. Notice that $\varphi$ is not an $[1,l]$-covering, otherwise we would have $\chi'_{[1,m]}(G) = \chi'_{[1,l]}(G)$, against our assumption.
Let $\varphi'$ be an equalized optimal $m$-bounded colouring of $G$, obtained using \lemma(DW). Necessarily $\varphi'$ is a $[t,t+1]$-covering of $G$, for some $t$. If $t<l$, then $\varphi'$ is a $[1,l]$-covering, which implies that $\chi'_{[1,m]}(G) = \chi'_{[1,l]}(G)$, a contradiction. Hence $t \geq l$. But then $\varphi'$ is an $[l,m]$-covering of $G$ and hence $\chi'_{[l,m]}(G)=\chi'_{[1,m]}(G)$ holds. Thus the algorithm in this case correctly outputs the excessive $[l,m]$-index and an excessive $[l,m]$-factorization of $G$.

Suppose now $\chi'_{[1,m]}(G) = \chi'_{[1,l]}(G)$:\\
if $\chi'_{[1,m]}(G) = \chi'_{[1,l]}(G)= \chi'_{[l]}(G)$, then the identities $\chi'_{[l,m]}(G)=\chi'_{[1,m]}(G)=\chi'_{[l]}(G)$ trivially hold by \theorem(XX).
In this case the algorithm correctly returns the excessive $[l,m]$-index $\chi'_{[l,m]}(G)$ and an excessive $[l,m]$-factorization of $G$ which is simply an excessive $[l]$-factorization (notice that $\chi'_{[l,m]}(G)$ is finite in this case since it equals $\chi'_{[1,m]}(G)$ which is always  finite).

Hence we can assume 
$$
k=\chi'_{[1,m]}(G) = \chi'_{[1,l]}(G) < \chi'_{[l]}(G).
$$
It will suffice to prove that no excessive $[1,m]$-factorization is an $[l,m]$-covering of $G$. This, by \theorem(XX), will imply that $\chi'_{[l,m]}(G)=\chi'_{[l]}(G)$ and hence prove the correctness of the algorithm.

Let $\cal M$ be an excessive $[1,m]$-factorization of $G$ and suppose that all matchings of $\cal M$ have cardinality at least $l$. In particular at least one matching of $\cal M$ has cardinality strictly larger than $l$, otherwise $\chi'_{[1,m]}(G)=\chi'_{[l]}(G)$, a contradiction.
Let $\cal L$ be an excessive $[1,l]$-factorization of $G$. At least one matching of $\cal L$ must have cardinality smaller than $l$, otherwise $\chi'_{[1,l]}(G) = \chi'_{[l]}(G)$, a contradiction.\\
Let $G_1$ and $G_2$ denote the multigraphs induced by $\cal M$ and $\cal L$, respectively. We have $$|E(G_2)|<lk$$ and hence $$|E(G)|<lk.$$
On the other hand $$|E(G_1)|>lk.$$
Thus, we can delete some edges from $G_1$ and still obtain a multigraph $H$ having exactly $lk$ edges and admitting $G$ as its underlying simple graph. Notice that $G_1$ (and hence $H$) is $k$-edge colourable by construction. An equalized $k$-edge-colouring of $H$ induces an excessive $[l]$-factorization of $G$. Hence $\chi'_{[m]}(G)=\chi'_{[l]}(G)$, a contradiction. 
This contradiction proves that no excessive $[1,m]$-factorization is an $[l,m]$-covering of $G$, and hence $\chi'_{[1,m]}(G)\neq \chi'_{[l,m]}(G)$.
By \theorem(XX) this implies that $\chi'_{[l,m]}(G)=\chi'_{[l]}(G)$ holds, and hence proves the correctness of our algorithm.\endproof


\end{document}